\documentclass[12pt,a4paper]{article}

\usepackage[utf8]{inputenc}
\usepackage[russian]{babel}
\usepackage{amssymb,amsmath}
\usepackage{graphicx}
\usepackage{color}
\usepackage{ulem}
\usepackage{hyperref}
\usepackage{extarrows}
\usepackage{floatflt}
\usepackage{caption}
\usepackage{mathtext}
\usepackage{amsfonts}
\usepackage{amsmath}
\usepackage{euscript}
\usepackage{amsthm}
\usepackage{framed}
\usepackage[justification=centering]{caption}

\usepackage{amssymb} 

\usepackage{geometry}

\usepackage{mathrsfs}

\newlength{\wdth}

\newif\ifhide
\hidefalse

\titlepage

\newtheorem{rTheorem}{Теорема}
\newtheorem{rlemma}{Лемма}
\newtheorem{rProposition}{Предложение}
\newtheorem{rdefinition}{Определение}
\newtheorem{rExample}{Пример}

\newtheorem{rCorollary}{Следствие}

\usepackage{hyperref}

\title{О центральной предельной теореме Добрушина для неоднородных цепей Маркова в схеме серий  \\
\author{
А.Ю.Веретенников\footnote{Институт проблем передачи информации РАН, Москва, Российская Федерация, alexander.veretennikov2011@yandex.ru}, А.И.Нуриева\footnote{Национальный исследовательский университет Высшая школа экономики, Москва, Российская Федерация, ai$\_$nurieva@mail.ru}}
}

\begin{document}

\maketitle

\begin{abstract}
В работе предложен новый вариант достаточного условия типа Добрушина в схеме серий для неоднородных цепей Маркова из достаточно широкого класса таких процессов. Само условие Добрушина может не выполняться, однако, ЦПТ при этом все же имеет место.

\medskip

\noindent
Ключевые слова: ЦПТ, перемешивание, схема серий.

\medskip

\noindent
Коды MSC: 60F05

\end{abstract}

\section{Введение}

В 1956 г. Р.Л. Добрушин установил важную центральную предельную теорему (далее ЦПТ) в схеме серий для цепей Маркова (далее ЦМ) с дискретным временем, вообще говоря, неоднородных  \cite{Dobr56}. Одним из главных предположений для такого результата было условие на скорость перемешивания, которая оценивается коэффициентом эргодичности, введенным в процессе выкладок еще самим А.А. Марковым. Позднее этот коэффициент получил имя: после публикации работы \cite{Dobr56} его стали называть коэффициентом эргодичности Добрушина. В попытке восстановить историческую справедливость в этой статье будем называть его коэффициентом Маркова -- Добрушина (далее МД). 
И до работы \cite{Dobr56} появлялись статьи на данную тему (см. \cite{Markov, Bern,  Sapogov, Linnik}), в основном, для важных, но все же частных классов цепей Маркова, причем зачинателем этой темы выступил сам А.А. Марков. Однако, только Добрушиным спустя почти полвека была установлена окончательная и в определенном смысле оптимальная асимптотика убывания коэффициента МД $\alpha_n$, достаточная для искомой ЦПТ в схеме серий, см. ниже условие (\ref{cor-e1}): 
$$
\lim_{n \to \infty} n^{1/3} \alpha_n = \infty.
$$
На примере, построенном еще ранее С.Н.Бернштейном \cite[глава 2]{Bern}, было показано, что если указанное условие на коэффициент МД не выполнено, то ЦПТ может не иметь места. В этом смысле условие Добрушина для общих ЦМ естественно считать окончательным.  Пример же этот в литературе часто называют примером Бернштейна -- Добрушина. 
Весьма непросто  написанная, статья Добрушина вызвала огромный интерес, прежде всего, благодаря основополагающей роли  ЦПТ в математической статистике, а также, очевидно, ввиду достижения наилучшей постоянной $1/3$ в достаточном условии и для общих цепей Маркова, в сочетании с упомянутым (контр)примером. Вскоре появилось много продолжений и обобщений этого результата, в частности, для сложных ЦМ, см. краткую заметку  \cite{Shi} о докладе на Московском математическом обществе. 
Уже не столь давно -- еще примерно через полвека после работы \cite{Dobr56} -- в статье 2005г. \cite{Seth-Varadhan}  основная теорема из \cite{Dobr56} была передоказана на основе существенно 
более простой мартингальной техники. 

В результате всех этих работ, -- и особенно ввиду того, что (контр)пример Бернштейна -- Добрушина построен на простейшем фазовом пространсте из двух точек, -- могло легко  возникнуть впечатление, что упомянутое выше (и далее в (\ref{cor-e1})) условие практически необходимо для ЦПТ, хотя логически такое заключение, конечно, было бы неверно. 
И вот, в связи с недавними результатами по обобщениям коффициента МД (см. \cite{VerBut, VerVer}), у одного из авторов  возникло желание применить эти обобщения и к задаче о ЦПТ в схеме серий. 
И действительно, в  данной работе будет показано, что условия работ \cite{Dobr56, Seth-Varadhan} в самом деле могут быть ослаблены за счет использования эргодического коэффициента не за один шаг цепи, а за два, или больше: как оказывается, столь простая  модификация может существенно расширить область  применимости предельной теоремы, разумеется, уже для более узкого класса ЦМ. В предпоследнем разделе будут приведены примеры ``непустоты'' класса моделей, удовлетворяющих новым условиям, для которого условия из \cite{Dobr56, {Seth-Varadhan}} не выполнены, однако, ЦПТ имеет место. При этом оказалось, что примеры удается построить лишь на фазовых пространствах с четырьмя и более элементами, так что возможно, что пример Бернштейна -- Добрушина в определенном смысле  неявно  препятствовал попыткам ослабить условия для ЦПТ для определенных классов ЦМ. 

Мотивация постановки вопроса следует из того же соображения, что  ЦПТ является одной из основ теоретической статистики, широко применяемой для оценки параметров, для построения доверительных множеств и для проверки гипотез. Стало быть, интересно любое расширение области ее применимости. 

Итак, будем рассматривать марковскую цепь с дискретным временем, не обязательно однородную. Пусть $\pi = \pi_{i,i+1} (x, dy)$ -- марковское переходное ядро на $(\mathbf{X}, \mathcal{B}(\mathbf{X}))$ в момент времени $i$, где $(\mathbf{X}, \mathcal{B}(\mathbf{X}))$ -- заданное измеримое пространство с обычными требованиями на марковский процесс (см. \cite{Dynkin}).
Под неоднородной марковской цепью длины $n$ на вероятностном пространстве $(\mathbf{X}, \mathcal{B}(\mathbf{X}))$ с переходной вероятностью $\{ \pi_{i,i+1} = \pi_{i,i+1} (x,dy): 1 \leq i \leq n-1 \}$ имеется в виду марковский процесс на произведении пространств $(\mathbf{X}^n, \mathcal{B}(\mathbf{X}^n))$ с вероятностями перехода  
$$
\mathsf P[X_{i+1} \in A | X_i = x]\stackrel{\text{п.н.}} = \pi_{i,i+1}(x,A),
$$
где $\{X_i: 1 \leq i \leq n \}$ -- траектория процесса на  отрезке времени $[1,n]$. В частности, при начальном\footnote{Отметим, что здесь удобно начальным моментом времени считать единицу, а не ноль.} распределении $X_1 \sim \mu$ распределение в момент $k \geq 2$ задается мерой $\mu\pi_{1,2}\pi_{2,3} \dots \pi_{k-1,k}$. При $i < j$ обозначим 
$$
\pi_{i,j} = \pi_{i,i+1}\pi_{i+1,i+2} \dots \pi_{j-1,j}.
$$

\begin{rdefinition}
Величина  $\delta(\pi)$, определенная одним из следующих двух эквивалентных выражений,
\begin{equation}\label{dpi}
\begin{aligned}
\delta(\pi) 
& = \sup_{\substack{x_1, x_2 \in \mathbf{X},\\ A \in \mathcal{B}(\mathbf{X})} }|\pi(x_1, A) - \pi(x_2, A)| 
 \\
&= \frac{1}{2} \sup_{\substack{x_1, x_2 \in \mathbf{X},\\ \|f\|_{L_\infty}\le 1}} | \int f(y) [\pi (x_1, dy) - \pi (x_2, dy)]|,
\end{aligned}
\end{equation}
называется коэффициент эргодичности Маркова -- Добрушина (в \cite{Seth-Varadhan} он назван просто коэффициентом сжатия). 
Положим также  
$$
\delta(\pi)(x_1,x_2)
= \sup_{A \in \mathcal{B}(\mathbf{X})} |\pi(x_1, A) - \pi(x_2, A)|  \\
= \frac{1}{2} \sup_{\|f\|_{L_\infty}\le 1)}| \int f(y) [\pi (x_1, dy) - \pi (x_2, dy)]|,
$$
$$
\alpha(\pi)(x_1,x_2):= 1- \delta(\pi)(x_1,x_2), \quad
\alpha(\pi) := 1 - \delta(\pi).
$$
Если понятно, какова переходная матрица, будем использовать сокращенные обозначения
$$
\alpha, \quad \delta, \quad
\alpha(x_1,x_2), \quad  \delta(x_1,x_2), 
$$
а для схемы серий -- также 
$$
\alpha_n, \quad \delta_n, \quad
\alpha_n(x_1,x_2), \quad  \delta_n(x_1,x_2).
$$
\end{rdefinition}

\noindent
Из определения видно, что $0 \leq \delta(\pi) \leq 1$ причем  $\delta(\pi) = 0$ тогда и только тогда, когда мера $\pi (x, dy)$ не зависит от $x$; в последнем случае вторая случайная величина не зависит от первой. Будем называть меру $\pi$ невырожденной, если $0 \leq \delta(\pi) < 1$. Для интегралов будем использовать  стандартные обозначения

$$
(\mu \pi)(A) = \int \mu(dx) \pi(x, A) \quad \text{и} \quad (\pi u)(x) \int  u(y) \pi(x, dy).
$$

\noindent
Отметим некоторые свойства коэффициента $\delta(\pi)$. Прежде всего, имеет место его эквивалентное представление, непосредственно вытекающее из первого определения в (\ref{dpi}), в виде  формулы
$$
\delta(\pi) = \sup_{\substack{x_1,x_2 \in \mathbf{X}, 
 \\ 
u \in \mathcal{U}}} |(\pi u)(x_1) - (\pi u)(x_2)|,
$$
где $\mathcal{U} = \{ u(\cdot): \sup_{y_1, y_2}|u(y_1) - u(y_2)| \leq 1\}$. В \cite{Seth-Varadhan} правая часть в этом представлении названа операторной нормой ядра $\pi$ по отношению к банаховой (квази-)норме $\operatorname{Osc}(u) = \sup_{x_1, x_2}|u(x_1) - u(x_2)|$, которую также называют осцилляторной (квази-)нормой функции $u$. Далее, для любых переходных мер $\pi_1, \pi_2$ имеет место неравенство 
$$
\delta(\pi_1 \pi_2) \leq \delta(\pi_1) \delta(\pi_2),
$$
где $\pi_1 \pi_2$ -- переходное ядро за два шага, когда первый совершается с помощью матрицы $\pi_1$, а второй -- с помощью $\pi_2$. Аналогично,
$$
\delta(\pi_{i,j}) \leq \delta(\pi_{i,i+1})\delta(\pi_{i+1,i+2}) \dots \delta(\pi_{j-1,j}).
$$
Соответственно, $\alpha(\pi_{i,j}) = 1 - \delta(\pi_{i,j})$.

Как обычно,  $\mathsf{E} Y$ и $\mathsf{D} [Y]$ обозначают, соответственно,   математическое ожидание и дисперсию случайной величины $Y$. 

В постановке схемы серий, для любого $n \geq 1$ через $\{ X_i^{(n)}: 1 \leq i \leq n \}$ обозначаем $n$ значений  неоднородной марковской цепи на множестве $\mathbf{X}$ с переходными  ядрами  $\{\pi^{(n)}_{i,i+1} = \pi^{(n)}_{i,i+1} (x, dy): 1 \leq i \leq n-1 \}$ и с начальным распределением $\mu^{(n)}$. Положим
$$
\alpha_n = \min_{1 \leq i \leq n -1} \alpha(\pi^{(n)}_{i,i+1}) , \quad
\alpha^{(2)}_n = \min_{1 \leq i \leq n - 2} \alpha(\pi^{(n)}_{i,i+2}).
$$
Далее,  пусть $\{ f^{(n)}_i: \mathbf{X} \mapsto \mathbb R, \, 1 \leq i \leq n \}$ -- вещественные функции на $\mathbf{X}$. Для каждого $n \geq 1$ положим
$$
S_n := \sum_{i=1}^n  f^{(n)}_i(X^{(n)}_i).
$$
Основной результат работы Добрушина вместе со следствием приведен в статье  \cite{Seth-Varadhan} в следующем несколько  упрощенном виде. 
\begin{rProposition}
Пусть для некоторой конечной постоянной $C_n$
$$
\sup_{1 \leq i \leq n} \sup_{x \in \mathbb{X}} |f^{(n)}_i(x)| \leq C_n.
$$
Тогда, если верно, что
\begin{equation}\label{thm-e1}
\lim_{n \to \infty} C^2_n \alpha_n^{-3} \big[ \sum_{i=1}^n \mathsf{D}(f^{(n)}_i(X^{(n)}_i)) \big]^{-1} = 0,
\end{equation}
то справедлива ЦПТ в схеме серий:
\begin{equation} \label{CLT0}
\frac{S_n - \mathsf{E}(S_n)}{\sqrt{{\mathsf D}(S_n)}} \implies {\cal N}(0,1).
\end{equation}
\end{rProposition}
Здесь ``$\implies$'' -- обозначение для слабой сходимости. 

\begin{rCorollary}\label{Cor1}
В частности, если  $\sup_{n} C_n = C < \infty$ и  $\inf_{i,n}\mathsf{D}(f^{(n)}_i(X^{(n)}_i)) \geq c > 0$, и если выполнено условие 
\begin{equation}\label{cor-e1}
\lim_{n \to \infty} n^{1/3} \alpha_n = \infty,
\end{equation}
то имеет место слабая сходимость (\ref{CLT0}). 
\end{rCorollary}
Итак, цель данной работы -- предложить ослабленные модификации условий (\ref{thm-e1}) и (\ref{cor-e1}), при которых ЦПТ (\ref{CLT0}) по-прежнему имеет место, хотя сами эти условия могут нарушаться. Работа не касается никаких уточнений к ЦПТ типа скорости сходимости ``на хвостах'', больших уклонений, и др., которым также посвящен целый ряд известных публикаций.

Статья состоит из шести разделов: 1 -- данное введение, 2 -- основные результаты, 3 -- вспомогательные утверждения, 4 -- доказательство теоремы \ref{thm1}, 5 -- примеры, 6 -- о примере Бернштейна -- Добрушина. 

\section{Основной результат}\label{sec2}
Прежде всего, заметим, что справедливы следующие элементарные неравенства: если $1\le i  \le j \le n$ и  значение $j-i$ является четным, то 
$$
\delta(\pi_{i,j}) \leq \delta(\pi_{i,i+2})\delta(\pi_{i+2,i+4}) \dots \delta(\pi_{j-2,j}) \leq (1 - \alpha^{(2)}_n)^{\frac{j-i}{2}},
$$
а если нечетным, то
$$
\delta(\pi_{i,j}) \leq \delta(\pi_{i,i+2})\delta(\pi_{i+2,i+4}) \dots \delta(\pi_{j-3,j-1})\delta(\pi_{j-1,j}) \le (1 - \alpha^{(2)}_n)^{\frac{j-1-i}{2}}.
$$
Следующая теорема является основным результатом данной работы. В разделе \ref{sec:examples} она будет дополнена некоторыми примерами.   
\begin{rTheorem}\label{thm1}
Пусть при всяком $n$ 
$$
\sup_{1 \leq i \leq n} \sup_{x \in \mathbb{X}} |f^{(n)}_i(x)| \leq C_n <\infty. 
$$
Тогда, если выполнено условие 
\begin{equation}\label{condition0}
\lim_{n \to \infty} C^2_n \alpha_n^{-1} (\alpha^{(2)}_n)^{-2} \big[ \sum_{i=1}^n \mathsf{D}(f^{(n)}_i(X^{(n)}_i)) \big]^{-1} = 0,
\end{equation}
то имеет место ЦПТ в схеме серий\footnote{Здесь обозначения таковы, что схема серий оказывается как бы скрытой, но это именно она, тогда как обозначения соответствуют таковым в \cite{Dobr56} и \cite{Seth-Varadhan}.}
\begin{equation} \label{CLT}
\frac{S_n - \mathsf{E}(S_n)}{\sqrt{\mathsf{D}(S_n)}} \implies {\cal N}(0,1).
\end{equation}
\end{rTheorem}

\begin{rCorollary}\label{Cor2}
Если  $\sup_{n} C_n = C < \infty$ и  $\inf_{i,n}\mathsf{D}(f^{(n)}_i(X^{(n)}_i)) \geq c > 0$, и если выполнено условие 
\begin{equation}\label{cor2-e1}
\lim_{n \to \infty} n \alpha_n (\alpha^{(2)}_n)^2 = \infty,
\end{equation}
то имеет место сходимость (\ref{CLT0}). 
\end{rCorollary}

\section{Вспомогательные утверждения}
В этом разделе приведем необходимые для доказательства теоремы вспомогательные утверждения. Кроме лемм \ref{lem1},  \ref{lem2} и \ref{lem4}, все они  взяты без изменений из работ \cite{Hall-Heyde} и \cite{Seth-Varadhan}; такие цитирутся без доказательств. Для простых же, но существенных для предлагаемого обобщения лемм  \ref{lem1}, \ref{lem2} и \ref{lem4} приведены доказательства, поскольку они содержат некоторые новые элементы. Обозначения взяты, по большей части, такие же, как в работе \cite{Seth-Varadhan}.
В частности, через $\|Z\|_{L^{\infty}}$ будет обозначаться супремум-норма случайной величины:
$$
\|Z\|_{L^{\infty}} := \sup_{\omega\in \Omega} |Z(\omega)|.
$$
Для неслучайной функции $f(x)$ супремум-норма обозначается как $\|f\|$.

\begin{rProposition}[см. \cite{Hall-Heyde}, Corollary 3.1]\label{pro2}
Пусть при каждом $n \geq 1$ процесс $\{(W^{(n)}_i: 0 \leq i \leq n \}$ является  мартингал относительно фильтрации $({\cal G}^{(n)}_i)$,  $W^{(n)}_0 = 0$, и пусть  $\xi^{(n)}_i := W^{(n)}_i - W^{(n)}_{i-1}$. Если выполнены условия 
\begin{equation}\label{clt_a}
\max_{1 \leq i \leq n} \|\xi^{(n)}_i\|_{L^\infty} \to 0 , \quad \sum_{i = 1}^n \mathsf{E}[(\xi^{(n)}_i)^2 | G^{(n)}_i] \to 1 \quad \text{в} \quad L^2
\end{equation}
(то есть, 
$\sup_{\omega}\max_{1 \leq i \leq n} |\xi^{(n)}_i (\omega)| \to 0$
и 
$\mathsf E\left|\sum_{i = 1}^n \mathsf{E}[(\xi^{(n)}_i)^2 | G^{(n)}_i] - 1\right|^2 \to 0$), 
то имеет место слабая сходимость
$$
W^{(n)}_n \implies {\cal N}(0,1).
$$
\end{rProposition}
Иные версии  классических ЦПТ для мартингалов, в том числе, при несколько более слабых условиях, можно найти в других трудах, например, см.  \cite[Теорема 5.5.8]{Lip-Shi}. 

~

Всюду в дальнейшем, следуя \cite{Seth-Varadhan}, будем для простоты -- то есть, чтобы не заботиться далее о центрировании -- предполагать, что функции $\{ f^{(n)}_i\}$ таковы, что 
$$
\mathsf{E}[f^{(n)}_i(X^{(n)}_k)] = 0, \quad 1 \leq i \leq n, \; n \geq 1.
$$
Положим
$$
Z^{(n)}_k = \sum_{i=k}^n \mathsf{E}[f^{(n)}_i(X^{(n)}_i)| X^{(n)}_k],
$$
так что 
$$
Z^{(n)}_k = 
\begin{cases}
f^{(n)}_i(X^{(n)}_k) + \sum_{i=k+1}^n \mathsf{E}[f^{(n)}_i(X^{(n)}_k)| X^{(n)}_k], \quad \text{при} \quad 1 \leq k \leq n-1, 
 \\ \\
f^{(n)}_n (X^{(n)}_n), \quad \text{при} \quad k = n.
\end{cases}
$$
Тогда для $1 \leq k \leq n-1$,

\begin{equation} \label{MartEq}
f^{(n)}_i(X^{(n)}_k) = Z^{(n)}_k - \mathsf{E}[Z^{(n)}_{k+1}| X^{(n)}_k].  
\end{equation}
Кроме того, для $1 \leq k \leq n-2$ это выражение можно представить в виде  $(Z^{(n)}_k - \mathsf{E}[Z^{(n)}_{k}| X^{(n)}_{k-1}]) + (\mathsf{E}[Z^{(n)}_{k}| X^{(n)}_{k-1}] - \mathsf{E}[Z^{(n)}_{k+1}| X^{(n)}_k])$. Стало быть, сумма $S_n$ допускает  следующее мартингальное представление:
$$
S_n = \sum_{i=1^n } f^{(n)}_i(X^{(n)}_i) = \sum_{k=2}^n [Z^{(n)}_k - \mathsf{E}[Z^{(n)}_{k}| X^{(n)}_{k-1}]] + Z^{(n)}_1.
$$
Как отмечается и в \cite{Seth-Varadhan}, это крайне полезное преобразование было предложено в статье \cite{Gordin}. 
В частности, в силу очевидной некоррелированности слагаемых, имеем, 
$$
\mathsf{D}(S_n) = \sum_{k=2}^n \mathsf{D}(Z^{(n)}_k - \mathsf{E}[Z^{(n)}_{k}| X^{(n)}_{k-1}]) + \mathsf{D}(Z^{(n)}_1).$$ 
Положим
\begin{equation}\label{XiEq}
    \xi^{(n)}_k = \frac{1}{\mathsf{D}(S_n)}[Z^{(n)}_k - \mathsf{E}[Z^{(n)}_{k}| X^{(n)}_{k-1}]].
\end{equation}
Тогда процесс $M^{(n)}_k = \sum_{l=1}^n  \xi^{(n)}_l$ является  мартингалом относительно фильтрации $\mathcal{F}^{(n)}_k = \sigma\{X^{(n)}_l :1 \leq l \leq k\}$ при $1 \leq k \leq n$. Доказательство ЦПТ, как в \cite{Seth-Varadhan}, так и данной работе  сводится к тому, чтобы аппроксимировать выражение $\frac{S_n}{\sqrt{\mathsf{D}(S_n)}}$ мартингалом $M^{(n)}_n$ и затем использовать ЦПТ для мартингал - разностей,  в данном случае, предложение  \ref{pro2}. Стало быть, задача состоит в том, чтобы   проверить оба условия в (\ref{clt_a}). Приведем необходимые вспомогательные утверждения. Далее, выражение вида $f^{(n)}_i\pi_{i,j}f^{(n)}_j$ означает не что иное, как произведение функций $f^{(n)}_i$ и $f^{(n)}_i\pi_{i,j}f^{(n)}_j$. В данном тексте выражения вида $\lfloor\frac{j-i}{2}\rfloor$  в степени обозначают целую часть числа $\frac{j-i}{2}$, тогда как скобки $[...]$  -- всегда ``обычные'' квадратные скобки.

Отметим, что для любых $i<j$, $n$, и для всяких $f^{(n)}_j$, справедливо неравенство

$$
\operatorname{Osc}(\pi_{i,j}(f^{(n)}_j)) \le \delta(\pi_{i,j})\operatorname{Osc}(f^{(n)}_j).
$$

\begin{rlemma}\label{lem1}
В условиях теоремы, 
для всех $1 \leq i \leq j \leq n $, если $(j-i)$ -- четное, то верно, что
$$
\|\pi_{i,j} f^{(n)}_j\|_{L^\infty} \leq 2C_n (1 - \alpha^{(2)}_n)^{\frac{j-i}{2}}, \quad \operatorname{Osc}(\pi_{i,j}(f^{(n)}_j)^2) \leq 2 C^2_n (1 - \alpha^{(2)}_n)^{\frac{j-i}{2}},
$$
а также, для $1 \leq l < i \leq j \leq n $,
\begin{enumerate}
\item[(a)] при четном $(i - l)$ имеем,
$$
\operatorname{Osc}(\pi_{l,i}(f^{(n)}_i\pi_{i,j}f^{(n)}_j)) \leq 6 C^2_n (1 - \alpha^{(2)}_n)^{\frac{i-l}{2}}(1 - \alpha^{(2)}_n)^{\frac{j-i}{2}};
$$
    
\item[(b)] при нечетном $(i-l)$ имеем,
$$
\operatorname{Osc}(\pi_{l,i}(f^{(n)}_i\pi_{i,j}f^{(n)}_j)) \leq 6 C^2_n (1 - \alpha_n)(1 - \alpha^{(2)}_n)^{\lfloor\frac{i-l}{2}\rfloor}(1 - \alpha^{(2)}_n)^{\frac{j-i}{2}}.
$$
\end{enumerate}

В случае если $(j-i)$ -- нечетное, справедливы неравенства 
$$
\|\pi_{i,j} f^{(n)}_j\|_{L^\infty} \!\leq\! 2C_n (1 - \alpha^{(2)}_n)^{\lfloor\frac{j-i}{2}\rfloor}(1 - \alpha_n), \; \operatorname{Osc}(\pi_{i,j}(f^{(n)}_j)^2) \!\leq\! 2 C^2_n (1 - \alpha^{(2)}_n)^{\lfloor\frac{j-i}{2}\rfloor}(1 - \alpha_n),
$$
а также, для $1 \leq i \leq j \leq n $:
\begin{enumerate}
\item[(c)] при четном $(i - l)$ имеем,
$$
\operatorname{Osc}(\pi_{l,i}(f^{(n)}_i\pi_{i,j}f^{(n)}_j)) \leq 6 C^2_n (1 - \alpha_n)(1 - \alpha^{(2)}_n)^{\frac{i-l}{2}}(1 - \alpha^{(2)}_n)^{\lfloor\frac{j-i}{2}\rfloor};
$$

\item[(d)] при нечетном же $(i-l)$ имеем, 
$$
\operatorname{Osc}(\pi_{l,i}(f^{(n)}_i\pi_{i,j}f^{(n)}_j)) \leq 6 C^2_n (1 - \alpha_n)^2(1 - \alpha^{(2)}_n)^{\lfloor\frac{i-l}{2}\rfloor}(1 - \alpha^{(2)}_n)^{\lfloor\frac{j-i}{2}\rfloor}.
$$

\end{enumerate}

\end{rlemma}

\begin{proof}
Рассмотрим случай  когда числа $j-i$  и $i-l$ -- оба четные. Так как $\|f^{(n)}_j\| \leq C_n$, то $\operatorname{Osc}(f^{(n)}_j) \leq 2C_n$. Из определений коэффициента МД $\delta$ и коэффициента  $\alpha^{(2)}_n$ имеем, 
$$
\operatorname{Osc}(\pi_{i,j}f^{(n)}_j) \leq \operatorname{Osc}(f^{(n)}_j) \delta(\pi_{i,j}) \leq 2C_n(1 - \alpha_n)^{\frac{j-i}{2}}.
$$
Поскольку  $\mathsf{E}[(\pi_{i,j}f^{(n)}_j)(X^{(n)}_j)]=\mathsf{E}[f^{(n)}_j(X^{(n)}_j)] = 0$, то первое неравенство леммы получается из оценок
$$
\|\pi_{i,j} f^{(n)}_j\|_{L^\infty} \leq \operatorname{Osc} (\pi_{i,j}f^{(n)}_j) \leq 2 C_n (1 - \alpha^{(2)}_n)^{\frac{j-i}{2}}.
$$
Второе неравенство леммы устанавливается аналогично. Для доказательства  третьего (случай (a)) оцениваем:
\begin{equation*}
\begin{aligned}
& \operatorname{Osc}(\pi_{l,i}(f^{(n)}_i\pi_{i,j}f^{(n)}_j)) \leq (1 - \alpha^{(2)}_n)^{\frac{i-l}{2}}  \operatorname{Osc}(f^{(n)}_i\pi_{i,j}f^{(n)}_j) 
 \\ \\
&\le        (1 - \alpha^{(2)}_n)^{\frac{i-l}{2}} \bigg(\operatorname{Osc}(f^{(n)}_i)\|\pi_{i,j} f^{(n)}_j\|_{L^\infty} + \|f^{(n)}_j\|_{L^\infty}  \operatorname{Osc}(\pi_{i,j}f^{(n)}_j)\bigg) 
 \\\\
&\le 6C^2_n (1 - \alpha^{(2)}_n)^{\frac{i-l}{2}} (1 - \alpha^{(2)}_n)^{\frac{j-i}{2}},
\end{aligned}
\end{equation*}
что и требовалось.

Остальные случаи четности и нечетности $j-i$ и $i-l$ рассматриваются аналогично. В нечетных случаях естественным образом появляется множитель $(1 - \alpha_n)$, в последнем случае в квадрате.
\end{proof}
Следующее утверждение дает нижнюю оценку для дисперсии суммы, играющую важную роль в доказательстве. Нам не удалось придумать ее улучшение через коэффициент $\alpha^{(2)}_n$, поэтому предложение приводится дословно из работы \cite{Seth-Varadhan} без доказательства; однако, и с коэффициентом  $\alpha_n$ оно все равно крайне полезно для дальнейшего.
\begin{rProposition}[\cite{Seth-Varadhan}, proposition 3.2]\label{pro3}
В условиях теоремы \ref{thm1}, 
для $n \geq 1$
$$
\mathsf{D}(S_n) \geq \frac{\alpha_n}{4} \sum^n_{i=1} \mathsf{D}(f^{(n)}_i(X^{(n}_i)).
$$
\end{rProposition}

\begin{rlemma}\label{lem2}
В условиях теоремы \ref{thm1} справедливо равенство 
\begin{equation}\label{lem2-e1}
\lim_{n \to \infty} \sup_{1 \leq k \leq n} \frac{\|Z^{(n)}_k\|_{L^{\infty}}}{\sqrt{\mathsf{D}(S_n)}} = 0.
\end{equation}
\end{rlemma}
\begin{proof}
Из предыдущей леммы имеем,
\begin{align*}
&\|Z^{(n)}_k\|_{L^{\infty}} \leq \sum^n_{i=k} \|\mathsf{E} [f^{(n)}_i(X^{(n)}_i)|X^{(n)}_k]\|_{L^{\infty}} 
 \\\\
&\leq 2C_n (1 - \alpha_n)^{\lfloor\frac{n-k}{2}\rfloor}\sum_{i=k}^n (1 - \alpha^{(2)}_n)^{\lfloor\frac{i-k}{2}\rfloor} 
\leq 4 C_n \frac{1}{1 - \sqrt{1 - \alpha^{(2)}_n}}.
\end{align*}
Последний переход следует из формулы для суммы геометрической прогрессии с множителем $0 \leq \sqrt{1 - \alpha^{(2)}_n} < 1$. 
Поясним, что ввиду использования функции целой части, фактически, под знаком суммы имеют место пары одинаковых слагаемых, так что в данном случае геометрическая прогрессия удваивается -- отсюда возникновение множителя $4$ в правой части итоговой оценки, вместо двойки. Множитель $(1 - \alpha_n)^{\lfloor\frac{n-k}{2}\rfloor}$ был оценен сверху единицей.

Поскольку при $0< \alpha \le 1$ всегда справедливо неравенство 
$\displaystyle
\frac{1}{1 - \sqrt{1 - \alpha}} \le \frac{2}{\alpha}
$, то заключаем, что  при достаточно малых $\alpha^{(2)}_n$ выполнено неравенство 
$$
\|Z^{(n)}_k\|_{L^{\infty}} \leq \sum^n_{i=k} \|\mathsf{E} [f^{(n)}_i(X^{(n)}_i)|X^{(n)}_k]\|_{L^{\infty}} \leq 8C_n \frac{1}{\alpha^{(2)}_n}.
$$
Применяя, наконец, оценку предложения \ref{pro3}, имеем,
$$
\sup_{1 \leq k \leq n} \frac{\|Z^{(n)}_k\|_{L^{\infty}}}{\sqrt{\mathsf{D}(S_n)}} \leq 32 \, C_n \Bigg( (\alpha^{(2)}_n)^2 \alpha_n \sum^n_{i=1} \mathsf{D}(f^{(n)}_i(X^{(n)}_i)) \Bigg)^{-1/2}.
$$
В силу условия (\ref{condition0}) отсюда следует (\ref{lem2-e1}),
что и требовалось.
\end{proof}

Следующая лемма  доказываются совершенно так же, как в работе  \cite[lemma 3.3]{Seth-Varadhan}, поэтому приводим ее без доказательства. 

\begin{rlemma} \label{lem3}
Пусть $\left\{Y_l^{(n)}: 1 \leq l \leq n\right\}$ и $\left\{\mathcal{G}_l^{(n)}: 1 \leq l \leq n\right\}$, для $n \geq 1$, соответственно последовательность неотрицательных случайных величин и $\sigma$-алгебр такие, что $\sigma\left\{Y_1^{(n)}, \ldots, Y_l^{(n)}\right\} \subset \mathcal{G}_l^{(n)}$. Предположим, что 
$$
\lim _{n \rightarrow \infty} \mathsf{E}\left[\sum_{l=1}^n Y_l^{(n)}\right]=1 \quad \text {и} \quad \sup _{1 \leq i \leq n}\left\|Y_i^{(n)}\right\|_{L^{\infty}} \leq \epsilon_n,
$$
где $\lim _{n \rightarrow \infty} \epsilon_n=0$. Также предположим
$$
\lim _{n \rightarrow \infty} \sup _{1 \leq l \leq n-1} \operatorname{Osc}\left(\mathsf{E}\left[\sum_{j=l+1}^n Y_j^{(n)} \mid \mathcal{G}_l^{(n)}\right]\right)=0.
$$
Тогда
$$
\lim _{n \rightarrow \infty} \sum_{l=1}^n Y_l^{(n)}=1 \quad \text \quad \mbox{в} \quad L^2.
$$
\end{rlemma}

Далее, положим 
$$
v_j^{(n)} :=  \mathsf{E}[(\xi_j^{(n)})^2 \mid \mathcal{F}_{j-1}^{(n)}].
$$ 
Эти случайные величины измеримы относительно сигма-алгебр $\mathcal{G}^{(n)}_j = \mathcal{F}^{(n)}_{j-1}$ для $2 \leq j \leq n $.

\begin{rlemma} \label{lem4}
В условиях теоремы \ref{thm1} имеет место сходимость 
$$
\sup _{2 \leq l \leq n-1}\, \sup_{\omega\in \Omega}\, \operatorname{Osc}\left(\mathsf{E}\left[\sum_{j=l+1}^n v_j^{(n)} \mid \mathcal{F}_{l-1}^{(n)}\right](\omega)\right)=o(1) .
$$
\end{rlemma}

\begin{proof}
 Из соотношения \ref{thm-e1} (из теоремы \ref{CLT0}) находим следующее. Поскольку в силу мартингального свойства $\mathsf E\left[\xi_r^{(n)} \xi_s^{(n)} \mid \mathcal{F}_u^{(n)}\right]=0$ для $r>s>u$ (из выражений \ref{MartEq} и \ref{XiEq}), то, стало быть, имеем,

\begin{align}\label{l4-e1}
&\mathsf{E}\left[\sum_{j=l+1}^n v_j^{(n)} \mid \mathcal{F}_{l-1}^{(n)}\right]= \mathsf{E} \left[\sum_{j=l+1}^n\left(\xi_j^{(n)}\right)^2 \mid \mathcal{F}_{l-1}^{(n)}\right] 
=  \mathsf{E}\left[\left(\sum_{j=l+1}^n \xi_j^{(n)}\right)^2 \mid X_{l-1}^{(n)}\right] 
 \nonumber \\ \nonumber \\ \nonumber 
&= \left(\mathsf{D}S_n\right)^{-1} \mathsf{E}\left[\left(\sum_{j=l+1}^n f_j^{(n)}\left(X_j^{(n)}\right)-\mathsf{E}\left[Z_{l+1}^{(n)} \mid X_l^{(n)}\right]\right)^2 \mid X_{l-1}^{(n)}\right] 
 \\ \nonumber \\ 
&=  \left(\mathsf{D}S_n\right)^{-1} \mathsf{E}\left[\left(\sum_{j=l+1}^n f_j^{(n)}\left(X_j^{(n)}\right)\right)^2 \mid X_{l-1}^{(n)}\right] 
-\left(\mathsf{D}S_n\right)^{-1} \mathsf{E}\left[\mathsf E\left[Z_{l+1}^{(n)} \mid X_l^{(n)}\right]^2 \mid X_{l-1}^{(n)}\right] .
\end{align}
В силу леммы \ref{lem2} последнее выражение -- вычитаемое в (\ref{l4-e1}) -- при $n\to\infty$ ограничено величиной $\sup _{2 \leq l \leq n-1} \left(\mathsf{D}S_n\right)^{-1} \left\|Z_{l+1}^{(n)}\right\|_{L^{\infty}}^2=o(1)$, и поэтому его колебание  также  $o(1)$ равномерно по  $\omega$. 

Чтобы оценить колебание первого слагаемого -- уменьшаемого в выражении (\ref{l4-e1}), -- запишем:

\begin{align*}
&\operatorname{Osc}\left(\left(\mathsf{D}S_n\right)^{-1} \mathsf{E}\left[\left(\sum_{j=l+1}^n f_j^{(n)}\left(X_j^{(n)}\right)\right)^2 \mid X_{l-1}^{(n)}\right]\right) \\ \\
&\leq \left(\mathsf{D}S_n\right)^{-1} \sum_{l+1 \leq j, m \leq n} \operatorname{Osc}\left(\mathsf{E}\left[f_j^{(n)}\left(X_j^{(n)}\right) f_m^{(n)}\left(X_m^{(n)}\right) \mid X_{l-1}^{(n)}\right]\right).
\end{align*}
Для $l+1 \leq j \leq m \leq n$ используем лемму \ref{lem1}: в случае четных значений $j - l + 1$ и $m - j$, получим (хотя в данном случае показатели степеней целые, все равно  будут  использованы  обозначения для целых частей $\lfloor...\rfloor$, поскольку такие же оценки получаются и в остальных случаях при нечетных $j - l + 1$ или $m - j$, и там уже без обозначения  целых частей не обойтись),
$$
\operatorname{Osc}\left(\mathsf{E}\left[f_j^{(n)}\left(X_j^{(n)}\right) f_m^{(n)}\left(X_m^{(n)}\right) \mid X_{l-1}^{(n)}\right]\right) \leq 6 C_n^2\left(1-\alpha^{(2)}_n\right)^{\lfloor\frac{j-l+1}{2}\rfloor}\left(1-\alpha^{(2)}_n\right)^{\lfloor\frac{m-j}{2}\rfloor} .
$$
Как уже сказано, в остальных  случаях оценки аналогичны, только  появляются один, или два множителя $(1 - \alpha_n)$, которые оцениваются сверху 
единицей. Таким образом, искомые колебания равномерно ограничены по $l$ выражением $\left(\mathsf{D}S_n\right)^{-1} C^2_n (\alpha^{(2)}_n)^{-2}$. Ввиду этого, а также в силу  предложения \ref{pro3} и условия \ref{condition0}, получаем утверждение леммы.
\end{proof}

\section{Доказательство теоремы \ref{thm1}}

\begin{proof}
Из леммы \ref{lem2} следует, что достаточно показать сходимость 
$$
M_n^{(n)} / \sqrt{\mathsf D \left(S_n\right)} \implies \mathrm{\cal N}(0,1).
$$
Указанное свойство следует из предложения \ref{pro2}), если показать, что 
\begin{enumerate}
\item $\sup _{2 \leq k \leq n}\left\|\xi_k^{(n)}\right\|_{L^{\infty}} \rightarrow 0$,
    
\item $\sum_{k=2}^n \mathsf{E}\left[\left(\xi_k^{(n)}\right)^2 \mid \mathcal{F}_{k-1}^{(n)}\right] \rightarrow 1 \quad \text{в} \quad L^2$.
       
\end{enumerate}
Здесь пункт (1) следует из леммы \ref{lem2}, а пункт (2) -- из леммы \ref{lem3} и пункта (1). Тем самым, теорема доказана. 
\end{proof}

\section{Примеры}\label{sec:examples}
Конечно, один из важнейших моментов предлагаемого обобщения  заключается в примерах. Покажем, что в самом деле существует достаточно широкий класс процессов, для которых коэффициент эргодичности за один шаг недостаточно мал, а вот за два шага -- уже достаточно для справедливости ЦПТ в схеме серий. 
Сперва просто продемонстрируем, что бывают  случаи, когда $\alpha_n=0$, а 
$\alpha^{(2)}_n>0$. 

\begin{rExample}
Рассмотрим матрицу переходной вероятности для однородной цепи Маркова
$$
\pi = 
\begin{pmatrix}
    1 - 2 \beta_n & \beta_n & \beta_n & 0 \\
    0 & \frac{1}{2} & \frac{1}{2} & 0 \\
    0 & 0 & \frac{1}{2} & \frac{1}{2} \\
    0 & 0 & 0 & 1
\end{pmatrix}
$$
с параметром $0< \beta_n < 1/2$.
\end{rExample}
Подсчет коэффициентов дает значения 
$$
\delta_n = \sup_{x_i, x_j \in X, A \in B(X)} |\pi(x_i, A) - \pi(x_j, A)| = \max(1 - 2 \beta_n, \frac{1}{2}, \frac{1}{2}, 1) = 1, \; \alpha_n = 1 - \delta_n = 0.
$$
Переходная матрица за два шага  выглядит следующим образом:
$$
\pi^{(2)} = 
\begin{pmatrix}
    (1 - 2 \beta_n)^2 & \beta_n(1 -  2 \beta_n) + \frac{1}{2} \beta& \beta_n(1 -  2 \beta_n) 
 + \beta_n & \frac{1}{2} \beta_n\\
    0 & \frac{1}{4} & \frac{1}{2} & \frac{1}{4} \\
    0 & 0 & \frac{1}{4} & \frac{3}{4} \\
    0 & 0 & 0 & 1
\end{pmatrix}
$$
Вновь элементарный подсчет показывает, что за два шага коэффициенты  таковы: $\delta^{(2)}_n = 1 - \frac{1}{2}\beta_n , \quad \alpha^{(2)}_n = \frac{1}{2}\beta_n >0$, тогда как $\alpha_n = 0$.

~

Этот пример подсказывает, как можно строить в большом количестве ``настоящие'' примеры. Возмутим какой-либо элемент данной матрицы $\pi$ малым возмущением, так, чтобы коэффициент $\alpha_n$ стал отличным от нуля, однако, оставался сколь угодно малым. Параметр же $\beta_n$ будем подбирать таким образом, чтобы имело место условие (\ref{cor2-e1}). Тогда теорема \ref{thm1} окажется применимой, хотя условия теоремы Добрушина могут быть не  выполнены.

~

\begin{rExample}\label{ex2}
Рассмотрим переходную матрицу на фазовом пространсте из четырех точек,
$$
\pi = 
\begin{pmatrix}
    1 - 2 \beta_n & \beta_n & \beta_n & 0 \\
    0 & \frac{1}{2} & \frac{1}{2} & 0 \\
    0 & 0 & \frac{1}{2} & \frac{1}{2} \\
    0 & \varepsilon_n & 0 & 1 - \varepsilon_n
\end{pmatrix}
$$
с параметрами 
$$
\beta_n = \frac{1}{n^{1/6}}, 
\quad 
\varepsilon_n = \frac{1}{n^{1/3}}.
$$
\end{rExample}
Для матрицы за два шага получим:
$$
\pi^{(2)} = 
\begin{pmatrix}
    1 - 2\beta_n + \beta_n^2 & \frac{3}{2}\beta_n - 2 \beta_n^2 & 2\beta_n - 2 \beta_n^2 &  \frac{1}{2} \beta_n \\
    0 & \frac{1}{4} & \frac{1}{2} & \frac{1}{4} \\
    0 & \frac{1}{2} \varepsilon_n & \frac{1}{4} & \frac{3}{4} + \frac{1}{2} \varepsilon_n\\
    0 & \frac{3}{2} \varepsilon_n - \varepsilon_n^2 & \frac{1}{2} \varepsilon_n & 1 - 2 \varepsilon_n + \varepsilon_n^2
\end{pmatrix}
$$
Как нетрудно убедиться прямым подсчетом, при таком выборе имеют место асимптотики  $\alpha_n \sim n^{-1/3}$, 
$\alpha^{(2)}_n \sim (1/2) n^{-1/6}$. Действительно, по одному из эквивалентных определений, 
$
\alpha_n = \min_{1\le i<j\le 4}\alpha_n(i,j) 
$, 
где $\alpha_n(i,j) = \sum_{k=1}^4 \pi(i,k)\wedge\pi(i,k)$. Имеем, 
\begin{align*}
&\alpha_n(1,2) = 2\beta_n/2 = \beta_n; \; 
 \alpha_n(1,3) = \beta_n/2; \; 
 \alpha_n(1,4) = \varepsilon_n\wedge\beta_n; \;
  \\\\
& \alpha_n(2,3) = 1/2\wedge 1/2 = 1/2; \;  
 \alpha_n(2,4) = \varepsilon_n/2; \;
\alpha_n(3,4) = (1-\varepsilon_n)/2.   
\end{align*}
Видим, что при $\varepsilon_n < \beta_n$ минимум достигается на паре $ \alpha_n(1,4)$ и равен $\varepsilon_n$. 

Аналогично для $\pi^2$ имеем при $\beta_n\to 0$ (считаем только аспимпотики), 
\begin{align*}
&\alpha^{(2)}_n(1,2) \!\sim \! 3\beta_n/2 + 
2\beta_n \!+\!  \beta_n/2 \!=\! 4\beta_n; \; 
 \alpha^{(2)}_n(1,3) \!\sim\! 2\beta_n \!+\! \beta_n/2 \!=\! \beta_n/2; \; 
 \alpha^{(2)}_n(1,4) \!\sim \!\beta_n; \;
  \\\\
& \alpha^{(2)}_n(2,3) \sim  1/4 + 1/4 = 1/2; \;  
 \alpha^{(2)}_n(2,4) \sim 1/4; \;
\alpha^{(2)}_n(3,4) \sim 3/4.
\end{align*}
Видим, что при $n\to \infty$ минимум достигается на паре $ \alpha^{(2)}_n(1,3)$ и асимптотически равен $\beta_n/2$. Стало быть, при указанном выборе $\beta_n$ и $\varepsilon_n$ получаем,  
$$
\frac{(\alpha_n)^{-1}(\alpha^{(2)}_n)^{-2}}{n} \sim \frac{(\varepsilon_n)^{-1}(\beta_n/2)^{-2}}{n} \sim 
\frac{4 n^{2/3}}{n} \to 0, 
$$
так что теорема \ref{thm1} (при выполнении еще условия (\ref{condition0})) применима. При этом 
$$
\frac{(\alpha_n)^{-3}}{n} \sim
\frac{n^{1}}{n} = 1, 
$$
так что условия теоремы Добрушина не выполнены.


\medskip

\begin{rExample}\label{ex3}
Рассмотрим переходную матрицу на фазовом пространсте из пяти точек,
$$
\pi = 
\begin{pmatrix}
1 - 2 \beta_n & \beta_n  & \beta_n & 0  & 0
 \\
0 & 1/3 & 1/3 & 0 & 1/3 
 \\
0 & 0 & 1/2 & 1/2  & 0 
 \\
\varepsilon_n & 0 & 0 & 1-\beta_n & \beta_n - \varepsilon_n  
 \\
\varepsilon_n & 0 & 0 & \beta_n - \varepsilon_n  & 1-\beta_n
\end{pmatrix}
$$
с параметрами 
$$
\beta_n = \frac{1}{n^{1/6}}, 
\quad 
\varepsilon_n = \frac{1}{n^{1/3}}.
$$
\end{rExample}
Тогда переходная матрица за два шага имеет вид
$$
\begin{pmatrix}
(1 \!-\! 2 \beta_n)^2 \!&\! \frac{4}{3} \beta_n \!-\! 2 \beta^2_n  \!&\! \frac{11}{6} \beta_n \!-\! 2 \beta^2_n \!&\! 0\!&\!\frac{5}{6} \beta_n  
 \\
\frac{\varepsilon_n}{3} \!&\! \frac{1}{9} \!&\! \frac{5}{18} \!&\! \frac{1}{3}(\beta_n \!-\! \varepsilon_n) \!&\! \frac{11}{18} \!-\! \frac{\beta_n}{3}
 \\
\frac{\varepsilon_n}{2} \!&\! 0 \!&\! \frac{1}{4} \!&\! \frac{1}{2}(\beta_n \!-\! \varepsilon_n)  \!&\! \frac{3}{4} \!-\!  \frac{\beta_n}{2}
 \\
\varepsilon_n (2 \!-\! 2\beta_n \!-\! \varepsilon_n)\!&\! \varepsilon_n \beta_n\!&\! \varepsilon_n \beta_n \!&\! (1\!-\!\beta_n)^2 \!+\!  (\beta_n \!-\! \varepsilon_n )^2 \!&\! 2(\beta_n \!-\! \varepsilon_n)(1 \!-\! \varepsilon_n)  
 \\
\varepsilon_n (2 \!-\! 2\beta_n \!-\! \varepsilon_n)\!& \!\varepsilon_n \beta_n \!&\! \varepsilon_n \beta_n \!&\! 2(\beta_n \!-\! \varepsilon_n)(1\! -\! \varepsilon_n) \!&\! (1\!-\!\beta_n)^2 \!+\!  (\beta_n \!-\! \varepsilon_n )^2 \!
\end{pmatrix}
$$
Подсчет показывает асимптотики, аналогичные тем, что в примере \ref{ex2}, только с немного иными константами:
$$
\alpha_n \sim \varepsilon_n =n^{-1/3}, \quad \alpha_n^{(2)} \sim 5\beta_n/6 = 5n^{-1/6}/6.
$$
Как и в примере \ref{ex2}, условие (\ref{cor2-e1}) выполнено, а (\ref{cor-e1}) -- нет. 

\medskip

\begin{rExample}\label{ex4}
Рассмотрим переходную матрицу 
$$
\pi = 
\begin{pmatrix}
1-\beta_n & \varepsilon_n  & \varepsilon_n  & \beta_n - 2\varepsilon_n 
 \\
\varepsilon_n  & 1-\beta_n & \beta_n - 2\varepsilon_n & \varepsilon_n 
 \\
 1/4 & 1/4 & 1/4 & 1/4 
 \\
 1/4 & 1/4 & 1/4 & 1/4 
\end{pmatrix}
$$
с параметрами 
$$
\beta_n = \frac{1}{n^{1/6}}, 
\quad 
\varepsilon_n = \frac{1}{n^{1/3}}.
$$
\end{rExample}
Здесь находим при $n\to\infty$ (элементы матрицы показаны с точностью до о-малых высшего порядка)
$$
\pi^2= 
\begin{pmatrix}
1-(7\beta_n/4) & \beta_n/4  & \beta_n/4  & 5\beta_n/4 
 \\
\beta_n/4  & 1-(7\beta_n/4) & 5\beta_n/4 & \beta_n/4
 \\
3/8 & 3/8 & 1/8 & 1/8 
 \\
3/8 & 3/8 & 1/8 & 1/8 
\end{pmatrix} .
$$
Подсчет показывает такие асимптотики (в первом случае просто равенство) для эргодических коэффициентов $\alpha_n$ и $\alpha^{(2)}_n$ при достаточно больших $n$: 
$$
\alpha_n = \alpha_n(1,2) = 4\varepsilon_n = 4n^{-1/3}, \quad 
\alpha^{(2)}_n = \alpha^{(2)}_n(1,2) \sim \beta_n = n^{-1/6}.
$$

{\bf Замечание.}
Ясно, что с той же целью дальнейшего ослабления условий можно использовать коэффициенты МД за три и более шагов ($\alpha^{(3)}_n$ и т.д.), во всяком случае, если условий с $\alpha_n$ и $\alpha^{(2)}_n$ окажется недостаточно. Соответствующие примеры, возможно,  потребуют матриц больших размерностей.

\section{О примере Бернштейна -- Добрушина}
Р.Л.Добрушин в работе \cite{Dobr56}  вслед за С.Н.Бернштейном \cite[
гл. 2, стр. 156 в русском издании Собр. соч.]{Bern} рассмотрел следующий пример, который иллюстрирует случай, когда $\alpha_n \sim n^{-1/3}$ и при этом ЦПТ в схеме серий не выполнена. Напомним этот пример вкратце, чтобы пояснить, почему для него наша теорема \ref{thm1} не дает никакого выигрыша при использовании коэффициента $\alpha^{(2)}_n$. Рассматривается сперва 
переходная матрица вероятностей с двумя состояниями $\{1,2\}$, 
$$
Q(p) = 
\begin{pmatrix}
    1 - p & p \\
    p & 1 - p
\end{pmatrix}\, .
$$
Чуть дальше величина $p$ будет сделана зависящей от времени и от номера серии, и затем определим $S_n = \sum_{i=1}^n 1(X_i=1)$ -- количество попаданий соответствующей цепи  Маркова в состояние $1$ за $n$ шагов. Определим также  $S_{a,b}= \sum_{i=a}^b 1(X_i=1)$. 

Пусть теперь величина $p=p_{i,n}$ принимает значения $\alpha_n$, либо $1-\alpha_n$, в зависимости от момента времени $i$, -- что  делает цепь неоднородной, -- и от номера серии $n$. Потребуем, чтобы $\alpha_n \to 0$ и при этом $\alpha_n \geq n^{-1/3}$. Положим  $k^{(n)}_i = i \lfloor(\alpha_n)^{-1}\rfloor$, $0\le i \le m_n$, где $m_n = \lfloor n/\lfloor\alpha_n^{-1}\rfloor\rfloor$, а также $m_{n+1} =n$. В примере Бернштейна -- Добрушина переходные матрицы в зависимости от момента времени задаются следующим образом:

$$
\pi^{(n)}_{i,i+1} = 
\begin{cases}
Q(\alpha_n), \quad \text{при} \quad i = 1 , \dots, k^{(n)}_1 - 1 \\
Q(1/2), \quad \text{при} \quad i = k^{(n)}_1 , \dots, k^{(n)}_{m_n}  \\
Q(1 - \alpha_n), \quad \text{для остальных} \quad i  \\
\end{cases}
$$
Матрицы $Q(1/2)$ в моменты времени $i = k^{(n)}_1 , \dots, k^{(n)}_{m_n}$ обеспечивают независимость частичных сумм на отрезках $[k^{(n)}_i, k^{(n)}_i + 1]$ при различных $i$. 

\medskip

Тогда, как оказывается (см. подробности в \cite{Dobr56, Seth-Varadhan}), имеют место асимптотики для дисперсий  $\mathsf{D}(S^{(n)}_{1, k^{(n)}_1}) \sim 4d_1(\alpha_n)^{-2}$ и  $\mathsf{D}(S^{(n)}_{k^{(n)}_1, n}) \sim d_2 n \alpha_n/2$, в обоих случаях с некоторыми положительными константами $d_1, d_2$. Из этого следует, что если $\alpha_n n^{1/3} \to \infty$, то $\alpha^{-2}_n \ll n \alpha_n$, и тогда  основной вклад в сумму $S^{(n)}$ вносит $S^{(n)}_{k^{(n)}_1, n}$, и так как это сумма большого числа независимых величин, то асимптотически (после полагающейся нормировки)  эта часть суммы распределена нормально. При этом сумма $S^{(n)}_{1, k^{(n)}_1-1}$ пренебрежимо мала. 

\medskip

Если же $\alpha_n = n^{-1/3}$, то $n \alpha_n = (\alpha_n)^{-2}$, и оказывается, что сумма $S^{(n)}_{1, k^{(n)}_1-1}$, не зависящая от  $S^{(n)}_{k^{(n)}_1, n}$, дает сравнимый --- не асимптотически пренебрежимый -- вклад в общую сумму. При этом предельное распределение у $S^{(n)}_{1, k^{(n)}_1-1}$ (опять же, после необходимой нормировки) другое -- невырожденное и  не гауссовское. В таком случае для всей суммы получается предельное распределение, являющееся композицией нормального и какого-то другого невырожденного  распределения. 
Как известно, такое распределение не может быть нормальным.

\medskip 
Почему же при $\alpha_n = n^{-1/3}$ не удается ослабить условие (\ref{cor-e1}), для получения итоговой ЦПТ, как это удалось в теореме \ref{thm1} и в примерах \ref{ex2} и \ref{ex3}? 
Оказывается, в примере Бернштейна -- Добрушина использование $\alpha^{(2)}_n$не может улучшить  ситуацию с ЦПТ  просто постольку, поскольку при $\alpha_n = n^{-1/3}$ асимптотики $\alpha^{(2)}_n$и $\alpha_n$ для матрицы $Q(\alpha_n)$ при малых $p_n$ имеют один и тот же порядок, что показывает простой подсчет. 
Стало быть, никакого выигрыша в асимптотике коэффициент $\alpha^{(2)}_n$ дать и не может. В определенном смысле, это свойство матриц $2\times 2$, а в больших размерностях оказывается вполне возможно, по крайней мере, для сильно разреженных матриц (см. предыдущий раздел), что при $\alpha_n \to 0$ выполнено  неравенство 
$\alpha^{(2)}_n \gg  \alpha_n$ (то есть, $\alpha_n/\alpha^{(2)}_n \to 0$), так что становится возможной ситуация когда $ n \alpha^3_n \not\to  \infty$, однако, $ n \alpha_n (\alpha^{(2)}_n)^2 \to  \infty$. В этой ситуации условия теоремы Добрушина не выполнены (см. следствие 1), тогда как теорема \ref{thm1} все равно применима. Во всяком случае, как демонстрируют примеры \ref{ex2}, \ref{ex3}, класс таких матриц достаточно велик и подобных примеров можно построить большое количество.

\section*
{Благодарности}
Для обоих авторов данная работа была поддержана грантом Фонда развития теоретической физики и математики «БАЗИС».

\end{document}